\documentclass[12pt]{article}
\usepackage{oldlfont}
\usepackage{enumerate}
\usepackage{amsmath}
\usepackage{amssymb}
\usepackage{amsthm}
\usepackage{mathrsfs}
\usepackage{amscd}
\usepackage{exscale}
\usepackage{latexsym}
\usepackage[all]{xy}
\usepackage{hyperref}
\usepackage{fancyhdr}
\usepackage{layout}
\usepackage{color}
\usepackage{graphicx}

\begin{document}

\baselineskip=17pt

\pagestyle{headings}

\numberwithin{equation}{section}

\makeatletter                                                           

\def\section{\@startsection {section}{1}{\z@}{-5.5ex plus -.5ex         
minus -.2ex}{1ex plus .2ex}{\large \bf}}                                 


\pagestyle{fancy}
\renewcommand{\sectionmark}[1]{\markboth{ #1}{ #1}}
\renewcommand{\subsectionmark}[1]{\markright{ #1}}
\fancyhf{} 
\fancyhead[LE,RO]{\slshape\thepage}
\fancyhead[LO]{\slshape\rightmark}
\fancyhead[RE]{\slshape\leftmark}

\addtolength{\headheight}{0.5pt} 
\renewcommand{\headrulewidth}{0pt} 

\newtheorem{thm}{Theorem}
\newtheorem{mainthm}[thm]{Main Theorem}

\newcommand{\ZZ}{{\mathbb Z}}
\newcommand{\GG}{{\mathbb G}}
\newcommand{\Z}{{\mathbb Z}}
\newcommand{\RR}{{\mathbb R}}
\newcommand{\NN}{{\mathbb N}}
\newcommand{\GF}{{\rm GF}}
\newcommand{\QQ}{{\mathbb Q}}
\newcommand{\CC}{{\mathbb C}}
\newcommand{\FF}{{\mathbb F}}

\newtheorem{lem}[thm]{Lemma}
\newtheorem{cor}[thm]{Corollary}
\newtheorem{pro}[thm]{Proposition}
\newcommand{\pf}{\noindent \textbf{Proof.} \ }
\newcommand{\eop}{$_{\Box}$  \relax}
\newtheorem{num}{equation}{}

\theoremstyle{definition}
\newtheorem{rem}[thm]{Remark}

\newcommand{\nsplit}{\cdot}
\newcommand{\G}{{\mathfrak g}}
\newcommand{\GL}{{\rm GL}}
\newcommand{\SL}{{\rm SL}}
\newcommand{\SP}{{\rm Sp}}
\newcommand{\LL}{{\rm L}}
\newcommand{\Ker}{{\rm Ker}}
\newcommand{\la}{\langle}
\newcommand{\ra}{\rangle}
\newcommand{\PSp}{{\rm PSp}}
\newcommand{\U}{{\rm U}}
\newcommand{\GU}{{\rm GU}}
\newcommand{\Aut}{{\rm Aut}}
\newcommand{\Alt}{{\rm Alt}}
\newcommand{\Sym}{{\rm Sym}}

\newcommand{\isom}{{\cong}}
\newcommand{\z}{{\zeta}}
\newcommand{\Gal}{{\rm Gal}}

\newcommand{\F}{{\mathbb F}}
\renewcommand{\O}{{\cal O}}
\newcommand{\Q}{{\mathbb Q}}
\newcommand{\R}{{\mathbb R}}
\newcommand{\N}{{\mathbb N}}
\newcommand{\E}{{\mathcal{E}}}

\newcommand{\DIM}{{\smallskip\noindent{\bf Proof.}\quad}}
\newcommand{\CVD}{\begin{flushright}$\square$\end{flushright}
\vskip 0.2cm\goodbreak}


\vskip 0.5cm

\title{On local-global divisibility by $p^2$\\ in elliptic curves}
\author{Laura Paladino, Gabriele Ranieri, Evelina Viada}
\date{  }
\maketitle

\vskip 1.5cm

\begin{abstract}
Let $p$ be a prime lager than $3$. Let $k$ be a number field, which does not contain the subfield of $\Q ( \zeta_{p^2} )$ of degree $p$ over $\Q$. Suppose that $\mathcal{E}$ is an elliptic curve defined over $k$. We  prove that the existence of a counterexample to the local-global divisibility by $p^2$ in $\mathcal{E}$, assures the existence of a $k$-rational point of exact order $p$ in $\mathcal{E}$. Using the Merel Theorem, we then shrunk the known  set of  primes for which there could be a counterexample to the local-global divisibility by $p^2$.
\end{abstract}

\section{Introduction}
Let $k$ be a number field  and let ${\mathcal{A}}$ be a commutative algebraic group defined over $k$.
Consider the following question:

\par\bigskip\noindent  P{\small ROBLEM}: \emph{Let $P\in {\mathcal{A}}(k)$. Assume that for all but finitely many  valuations $v\in k$, there exists $D_v\in {\mathcal{A}}(k_v)$ such that $P=qD_v$, where $q$ is a positive integer. Is it possible to conclude that there exists $D\in {\mathcal{A}}(k)$ such that $P=qD$?}

\par\bigskip\noindent This problem is known as \emph{Local-Global Divisibility Problem}.  By  B\'{e}zout's identity,  to get answers for a general integer $q$ it is sufficient to solve it for $q=p^n$, with $p$  a prime. In the classical case of ${\mathcal{A}}={\mathbb{G}}_m$, the answer is positive for $q$ odd, and negative for instance for $q=8$ (and $P=16$) (see for example \cite{AT}, \cite{Tro}).

\bigskip During the last few years R. Dvornicich and U. Zannier  \cite{DZ}  have given some general vanishing cohomological conditions, sufficient for a positive answer to the question; this led to a number of examples and counterexamples to the local-global principle  when  ${\mathcal{A}}$ is an elliptic curve or a torus of higher dimension (for tori see also  \cite{Ill} for further examples). \\

We are particularly interested in
 the case of  an elliptic curve $\mathcal{E}$. Dvornicich and Zannier \cite{DZ}  prove that the local-global principle holds on   $\mathcal{E}$ for  any prime (see also \cite{Won}).
Moreover,   if $\mathcal{E}$ does not admit a $k$-rational isogeny of degree $p$,
 then they give an affirmative answer when $q=p^n$, for all $n\in \mathbb{N}$ (see \cite{DZ3}).  By  J.-P. Serre \cite{Ser},  such an isogeny exists only for
 $p < C_{{\rm serre}}(k,{\mathcal{E}})$, with  $C_{{\rm serre}}(k,{\mathcal{E}})$ a constant  depending on $k$ and $\mathcal{E}$.  For  $\mathcal{E}$ defined over $\Q$, B. Mazur \cite{Maz1} proved that a rational isogeny of degree $ p $ exists only for  $p\in S=\{2,3,5,7,11,13,17,19,37,43,67,163\}$.
 Thus, by the mentioned results of   Dvornicich and  Zannier,
  in elliptic curves, the local-global divisibility by $p^n$ holds for all $n\in \mathbb{N}$ and $p>C(k,{\mathcal{E}})$;
    in elliptic curves over $\QQ$ it sufficies $p\notin S$. An open question is  if $S$ is minimal for the local-global question. Only for  $2^n$ and  $3^n$, for all $n\geq 2$, there are counterexamples, some explicit (see
\cite{DZ2}, \cite{Pal}, \cite{Pal2} and \cite{Pal3}).\\

    In  \cite{DZ3} page 30, Dvornicich and  Zannier  stated  the following question:\\
      {\it ``On the other hand, we do not know whether one may take  $C(k,{\mathcal{E}})$ independently of the curve $\E$. An analogue of the quoted result of Mazur would be helpful. However, to our knowledge, even the deep result of Merel \cite{Mer} do not cover the general case of rational cyclic groups, but only the case of rational points".}

In this paper we  prove such independence for local-global divisibility by $p^2$. We show that a counterexample to the local-global problem for $p^2$ implies the existence of a $k$-rational  point of exact order $p$, and not only of  a $k$-isogeny of degree $p$. As a consequence, we show that the set $S$ is not minimal.
\begin{thm}\label{main}

Let $p>3$ be a prime  and $\zeta_{p^2} $ a primitive $p^2$th root of unity. Let $\E$ be an elliptic curve defined over a number field $k$, which does not contain the subfield of $\Q ( \zeta_{p^2} )$ of degree $p$ over $\Q$.
Suppose that $\E$ does not admit any   $k$-rational torsion point of exact order $p$. Then,
if $P \in \E ( k )$ is a point, which is locally divisible by $p^2$ in $\E ( k_v )$ for all but finitely many valuations $v$,  $P$ is divisible by $p^2$ in $\E ( k )$.
\end{thm}

Then the quoted result of L. Merel in general,  the famous effective Mazur's theorem  \cite{Maz2}  for $\QQ$ and  the theorem of S. Kamienny \cite{Ke-Mo}  for quadratic number fields, produce positive and respectively effective answer  to the  question of Dvornicich and  Zannier for divisibility by $p^2$. In particular,  when we consider the local-global divisibility by $p^2$ in elliptic curves defined over $\Q$,
the set $S$ is not minimal and it can be shrunk to $\widetilde{S}=\{2,3,5,7\}$.

\begin{cor}\label{cor:cor42}
 Let $\E$ be an elliptic curve defined over any number field $k$.
Then, there exists a constant ${C} ( k )$, depending only on $k$, such that, for every prime number $p > {C} ( k )$, if $P \in \E ( k )$ is a point, which is locally divisible by $p^2$ in $\E ( k_v )$ for all but finitely many valuations $v$, then $P$ is divisible by $p^2$ in $\E ( k )$. In particular, ${C} ( \QQ )=7$ and  for a quadratic number field ${C} ( k )=13$.
\end{cor}

\DIM  By  \cite{Mer}, for every number field $k$, there exists a constant $C_{{\rm merel}} ( k )$ depending only on $k$, such that, for every prime $p > C_{{\rm merel}} ( k )$, no elliptic curve defined over $k$  has a $k$-rational torsion point of exact order $p$.
Let $p_0$ be the largest prime such that $k$ contains the subfield of $\Q ( \zeta_{p_0^2} )$ with degree $p_0$ over $\Q$. Observe that $p_0 \leq [k: \Q]$.
Set
\[
{C} ( k ) = {\rm max} \{ 3, p_0, C_{{\rm merel}} ( k ) \}.
\]
\noindent Then, apply Theorem~\ref{main}.

By  \cite{Maz2}, no elliptic curve defined over $\Q$ has a rational point of exact prime order larger than $7$. By~\cite{Ke-Mo}, no elliptic curve defined over a quadratic number field $k$  has a $k$-rational point of exact prime order larger than $13$.
\CVD

 In order to prove Theorem \ref{main}, we use several of the results of  Dvornicich and Zannier and some of the classical results of Serre (see Section~\ref{sec2}). The existence of a counterexample produces a $k$-isogeny of degree $p$, which in turn gives one unique or two distinct $\Gal ( k(\E[p]) / k )$-submodules of $ \E[p]$. The resulting simplicity of  the group $\Gal ( k(\E[p]) / k )$ makes it possible to generate a $k$-torsion point of exact order $p$. The proof of this main result is presented in Section~\ref{sec3}.

\bigskip\noindent \emph{Acknowledgments}. We would like to kindly  thank A. Bandini and F. Gillibert for their remarks and helpful discussions. The second and third author thank the FNS for financial support.
We  are grateful to the University of Basel and to the University of Calabria for the hospitality.

\section{Preliminary results}\label{sec2}

Let $k$ be a number field and let $\mathcal{E}$ be an elliptic curve defined over $k$. Let $p $ be a prime not $2$ or $3$. For  every positive integer $n$, we denote by ${\mathcal{E}}[p^n]$  the $p^n$-torsion subgroup of ${\mathcal{E}}$ and by $K_n = k({\mathcal{E}}[p^n])$  the number field obtained by adding to $k$ the coordinates of the $p^n$-torsion points of ${\mathcal{E}}$.
Let $G_n = \Gal ( K_n / k )$. As usual, we shall view   $\E[p^n]$ as $ \Z/p^n\Z \times \Z/p^n\Z$ and consequently we shall represent $G_n$ as a subgroup of
 $\textrm{GL}_2 ( \Z / p^n \Z )$, denoted by the same symbol.  By Silverman \cite[Chapter III, Corollary 8.1.1]{Sil}, the field  $ K_n$ contains a primitive $p^n$th root of unity $\zeta_{p^n} $. The basic properties of the Weil pairing  \cite[ III.8]{Sil} entail that the action of $G_n$ on $\zeta_p$ is given by $g(\zeta_p)=\zeta_p^{\det g}$ for $g \in G_n$.

As mentioned, the answer to the \emph{Local-Global Divisibility Problem} when $q=p^n$ is strictly connected to the vanishing condition of the  cohomological  group $H^1(G_n, {\mathcal{E}}[p^n])$ and of the local cohomological group $H^1_{\rm loc}(G_n,{\mathcal{E}}[p^n])$.  Let us recall some definitions and results for $\E$.
\par\bigskip

\noindent \textbf{Definition}[Dvornicich, Zannier \cite{DZ}] \label{defloc}\emph{Let $\Sigma$ be a group and let $M$ be a $\Sigma$-module.
We say that a cocycle $[c]=[\{Z_{\sigma}\}]\in H^1(\Sigma,M)$ satisfies the \emph{local conditions} if there exists $W_{\sigma}\in M$ such that $Z_{\sigma}=(\sigma - 1)W_{\sigma}$, for all $\sigma\in \Sigma$.
We denote by $H^1_{{\rm loc}}(\Sigma,M)$ the subgroup of $H^1(\Sigma,M)$ formed by such cocycles.
Equivalently, $H^1_{{\rm loc}} ( \Sigma, M )$ is the intersection of the kernels of the restriction maps $H^1 ( \Sigma, M ) \rightarrow H^1 ( C, M )$ as $C$ varies over all cyclic subgroups of $\Sigma$.}

\bigskip\noindent Working with all valuations, instead of almost all, one would get the classical definition of the Shafarevich group. 

\begin{thm}[Dvornicich, Zannier \cite{DZ}] \label{si1}\par \noindent Assume that
$H^1_{\rm loc}(G_n,{\mathcal{E}}[p^n])=0.$
 Let $P\in {\mathcal{E}}(k)$ be a point locally divisible by $p^n$ almost everywhere in the completions $k_v$ of $k$. Then there exists a point $D\in {\mathcal{E}}(k)$, such that $P=p^nD$.
\end{thm}

In \cite{DZ3} they prove that this theorem is not invertible.  Moreover, an intrinsic version of their main theorem  \cite{DZ3} (remark just after the main theorem and the first few lines of its proof) gives:
\begin{thm}\label{teo:teo11}  Suppose that $\E$ does not admit any $k$-rational isogeny of degree $p$.
Then $H^1 (G_n, \E[p^n] ) = 0$, for every $n \in \N$.
\end{thm}
 The next lemma is essentially  proved  in the proof of the Theorem \ref{teo:teo11}, in  \cite{DZ3} beginning of page 29.
 \begin{lem}\label{lem:lem12}
Suppose that there exists  a nontrivial multiple of the identity $\tau \in G_1$.
Then $H^1 ( G_n, \E[p^n] ) = 0$, for every $n \in \N$.
\end{lem}

\noindent  Another remark along their proof, concerns the group $G_D = G_1 \cap D$, where $D$  is the subgroup  of  diagonal matrices of $\textrm{GL}_2 ( \F_p )$.

\begin{cor}\label{cor:cor13}
Suppose that $G_D$ is not cyclic. Then $H^1 ( G_n, \E[p^n] ) = 0$, for every $n\in \N$.

\end{cor}
\DIM Since $G_D$ is not cyclic and
$
D \cong \Z/ ( p-1 ) \Z \times \Z/ ( p-1 ) \Z,
$
then $ G_D $ contains all the diagonal matrices of order $l$, for $l$ dividing $p-1$. Apply Lemma~\ref{lem:lem12}.\CVD

We now  describe $G_1$, under the hypothesis  $H^1 ( G_2, \E[p^2] ) \neq 0$.

\begin{lem}\label{lem:lem14}
Suppose  that $H^1 ( G_2, \E[p^2] ) \neq 0$. Then:
\begin{itemize}
\item[(i)] Either   $ \E[p] $ contains two distinct $G_1$-modules  of order $ p $ and  $G_1=\langle \rho \rangle$ where $\rho$ is a diagonal matrix, in a suitable $\F_p$-basis  of $\E[p]$.
\item[(ii)] Or  $\E[p]$ contains a unique $G_1$-module  of order $ p $
 and $ G_1 =\langle \rho,  \sigma \rangle$ with  $\rho$ a diagonal matrix and $  \sigma={1\,\,\, 1 \choose 0\,\,\, 1}$, in a suitable  $\F_p$-basis of $\E[p]$. Furthermore, for  $k^\prime \subseteq K_1$  the field fixed by $\rho$,  $[k^\prime : k] = p$ and  $k^\prime / k$ is Galois if and only if $\rho$ is the identity.
\end{itemize}

\end{lem}

\DIM

By Theorem~\ref{teo:teo11}, $\E$ admits a $k$-rational isogeny $\phi$ of degree $p$.
Since $\ker ( \phi )$ is a $G_1$-module, then $\E[p]$ contains either  one unique or  two distinct $G_1$-submodules of order $p$.

Part (i)

Let  $C_1=<P_1>$ and $ C_2=<P_2>$  be the two distinct cyclic $G_1$-submodules of order $p$.
Then, for every $\tau \in G_1$, there exist $\lambda_\tau, \mu_\tau$ such that
\[
\tau ( P_1 ) = \lambda_\tau P_1, \ \tau ( P_2 ) = \mu_\tau P_2.
\]
Hence, in the basis $\{P_1, P_2 \}$, $G_1$ is contained in the group of  diagonal matrices of $\GL_2 (\Z/p\Z)$.
By Corollary~\ref{cor:cor13},  $G_1$ is cyclic.\\

Part (ii)

 Let $P_1$ be a generator of the unique $G_1$-module  $C$ of order $ p $. Then  for every $\tau \in G_1$, there exists $\lambda_\tau \in \F_p$ such that $P_1^\tau = \lambda_\tau P_1$. Thus, for any  basis $\{P_1, P_2\}$, every element of $G_1$ is an upper triangular matrix.

We first show that $\sigma \in G_1$.
Since $C$ is unique, $G_1$ is not  contained in the group $D$ of the diagonal matrices.
Therefore, there exists
\[
\gamma =
\left(
\begin{array}{cc}
a & b \\
0 & d \\
\end{array}
\right)
\]
in $G_1$, with $a, b, d \in \F_p$ and $b \neq 0$.

If $a = d$, then
\[
\gamma^{p-1} =
\left(
\begin{array}{cc}
1 & \hspace{0.2cm} ( p - 1 ) a^{p-2} b \\
0 & 1 \\
\end{array}
\right),
\]
which implies that $\sigma \in \langle \gamma \rangle$.

If $a \neq d$, then $\gamma$ is diagonalizable. In   a basis of eigenvectors $\{P_1,P'_2\}$
\[
\gamma =
\left(
\begin{array}{cc}
a & 0 \\
0 & d \\
\end{array}
\right).
\]
Since $G_1 \not \subseteq D$, there exists $\delta \in G_1$ such that,  in the basis of eigenvectors $\{P_1,P'_2\}$,
\[
\delta =
\left(
\begin{array}{cc}
a^\prime & b^\prime \\
0 & d^\prime \\
\end{array}
\right),
\]
where $a^\prime, b^\prime, d^\prime$ are in $\F_p$ and $b^\prime$ is non-zero.
A simple computation shows that
\[
\delta \gamma \delta^{-1} \gamma^{-1} =
\left(
\begin{array}{cc}
1 & ( d - a ) b^\prime / d d^\prime \\
0 & 1 \\
\end{array}
\right).
\]
Since $b^\prime \neq 0$ and $a \neq d$, we get $\sigma \in \langle \delta \gamma \delta^{-1} \gamma^{-1} \rangle$.
So we have proved that $\sigma \in G_1$. Now we look for $\rho$.
Recall that $G_D = G_1 \cap D$.  By Corollary~\ref{cor:cor13}, $G_D=\langle \rho \rangle$ with $\rho$ diagonal and of order dividing $p-1$.  Let
\[
\sigma_1 =
\left(
\begin{array}{cc}
x & y \\
0 & z \\
\end{array}
\right)\in G_1.
\]
Observe that
\[
\sigma_2 =
\left(
\begin{array}{cc}
1 & -y / x \\
0 & 1 \\
\end{array}
\right)
\in \langle \sigma \rangle
\]
and
\[
\sigma_1 \sigma_2 =
\left(
\begin{array}{cc}
x & 0 \\
0 & z \\
\end{array}
\right).
\]
Then, every right coset of $G_D$ in $G_1$, contains an element of $\langle \sigma \rangle$.
Moreover by a simple verification,  such an element is unique.
Then
\begin{equation}\label{eqn:rel11}
\vert G_1 \vert / \vert G_D \vert = \vert \langle \sigma \rangle \vert = p.
\end{equation}
Since  $\sigma$ has order $p$ and $\rho$ has order   $\vert G_D \vert$ coprime to $p$, then $\rho$ and $\sigma$ generate $G_1$, which has order $p \vert G_D \vert $ .

By definition $k^\prime$ is fixed by $G_D=\langle \rho \rangle$,  thus  relation~(\ref{eqn:rel11}) gives $[k^\prime : k] = p$.

By Lemma~\ref{lem:lem12}   either $\rho$ is  the identity or  $\rho \neq \lambda I$, with $\lambda \in \F_p^\ast$.
If $\rho$ is the identity, we have $K_1 = k^\prime$ and $k^\prime / k$ is Galois.
If $\rho$ is not a scalar multiple of the identity, it is an easy computation to verify that $\sigma \rho \sigma^{-1}$ is not diagonal.
Thus $G_D$ is not normal in $G_1$.
Then, by Galois correspondence,  $k^\prime$ is not Galois over $k$.
  \CVD

Let $H = \Gal ( K_2 / K_1 )$.
Then $H$ fixes $\E[p]$
and  $\tau \in H$ is congruent to the identity modulo $p$.
Hence $H$ is contained in the subgroup of $\textrm{GL}_2 ( \Z / p^2 \Z )$, given by the elements
\[
\left(
\begin{array}{cc}
1 & 0 \\
0 & 1 \\
\end{array}
\right)
+
p
\left(
\begin{array}{cc}
a & b \\
c & d \\
\end{array}
\right),
\]
where $a, b, c, d \in \Z / p \Z$.
Such a subgroup is a $\F_p$-vector space of dimension $4$.
Then $H$ is a $\F_p$-vector space of dimension $\leq 4$.

In  \cite[Proposition 3.2]{DZ} and \cite[Proposition 2.5]{DZ} it is shown that  what really matters for the triviality of the first cohomological group is the $p$-Sylow. Since the structure of our group $G_1$ is particularly simple, we can apply their propositions to deduce:
\begin{lem}\label{lem:lem16}
If ${\rm dim}_{\F_p} ( H ) \neq 2$, then $H^1_{{\rm loc}} ( G_2, \E[p^2] ) = 0$.
\end{lem}

\DIM Assume that $H^1_{{\rm loc}} ( G_2, \E[p^2] ) \not= 0$, then $H^1 ( G_2, \E[p^2] ) \neq 0$.
By Lemma~\ref{lem:lem14}, either $G_1 \cong G_2/ H$ is generated by a diagonal matrix, or $G_1$ is generated by a diagonal matrix and
\[
\sigma =
\left(
\begin{array}{cc}
1 & 1 \\
0 & 1 \\
\end{array}
\right).
\]
In both cases $ H $ is contained in the $p$-Sylow  $H_p$ of $G_2$ (which is $H$ in the first case, and it is generated by $H$ and any extension of $\sigma$ to $K_2$ in the second case).
Then $H$ is the intersection of  $H_p$ and the kernel of the reduction modulo $p$ from $\textrm{G}L_2 ( \Z / p^2 \Z )$ to $\textrm{GL}_2 ( \Z / p \Z )$.
 In~\cite[Proposition 3.2]{DZ}, it is proven that: if  $p \neq 2, 3$ and ${\rm dim} ( H ) \neq 2$,  then $H^1_{{\rm loc}} ( H_p, ( \Z / p^2 \Z )^2 ) = 0$. Identifying $\E[p^2] $ with $ (\Z / p^2 \Z )^2$ we have $H^1_{{\rm loc}} ( H_p, \E[p^2] ) = 0$.
In~\cite[Proposition 2.5]{DZ}, it is shown that the triviality of the local cohomology group on the p-Sylow implies $H^1_{{\rm loc}} ( G_2, \E[p^2] ) = 0$. This shows our lemma.
\CVD

\section{Local-global divisibility by $p^2$} \label{sec3}

In this section we prove the central propositions of our work, these are the heart of the proof of the main theorem, presented at the end of the section.
First of all,  we recall a well known result that we shall  use in the following.

\begin{lem}\label{lem:lem21}
Let $\Delta$ be a cyclic group and let $M$ be a $\Delta$-module.
Set $N_{\Delta} \colon M \rightarrow M$
the homomorphism that sends $m \in M$ into $\sum_{i = 0}^{\vert \Delta \vert - 1} \delta^i m$,
where $\delta$ is a generator of $\Delta$.
Let
\[
{\rm Im} ( \delta -1 ) = \{m \in M \ \mid \ \exists \ m^\prime \in M \ \textrm{such} \ \textrm{that} \ m = \delta ( m^\prime ) - m^\prime\}.
\]
Then $H^1 ( \Delta, M ) \cong \ker ( N_{\Delta} ) / {\rm Im} ( \delta - 1 ).$
\end{lem}

\DIM See \cite[Chapter II, Example 1.20]{Mil}.
\CVD

\begin{pro}[The cyclic case]\label{pro:pro22} Suppose that $k$  does not contain the subfield  $L$ of $\Q ( \zeta_{p^2} )$ such that $[L : \Q] = p$.
Suppose  that $H^1_{{\rm loc}} ( G_2, \E[p^2] ) \neq 0$ and that $G_1$ is cyclic.
Then $K_1 = k ( \zeta_p )$, and there exists a $k$-rational point of $\E$ of exact order $p$.
\end{pro}

\DIM Recall that $G_1$ has order $d$ dividing $p-1$.
Consider the group $\Gal ( K_1 ( \zeta_{p^2} ) / k )$.
Since $k$ does not contain $L$, then  $kL / k$ is a cyclic extension of degree $p$.
Moreover $k L \cap K_1 = k$, because $[kL: k]$ and $[K_1: k]$ are coprime.
Finally $K_1 ( \zeta_{p^2} ) = K_1 L$, because $\zeta_p \in K_1$.
Then, by elementary Galois theory, we have that
\begin{equation}\label{eqn:rel21}
\Gal ( K_1 ( \zeta_{p^2} ) / k ) \cong \Gal ( kL / k ) \times \Gal ( K_1 / k ) \cong \Z / p \Z \times \Z / d \Z \cong \Z / pd \Z
\end{equation}
is a cyclic group.
Recall that $H = \Gal ( K_2 / K_1 )$. By Lemma~\ref{lem:lem16},  $\vert H \vert = p^2$. In addition, we know that $\zeta_{p^2} \in K_2$. Then $K_2 / K_1 ( \zeta_{p^2} )$ is a cyclic extension of degree $p$.
Let $\gamma$ be a generator of $\Gal ( K_2 / K_1 ( \zeta_{p^2} ) )$.
We obtain the following inflation-restriction sequence:
\begin{equation}\label{eqn:rel22}
0 \rightarrow H^1 ( \Gal ( K_1 ( \zeta_{p^2} ) / k ), \E[p^2]^{\langle \gamma \rangle} ) \rightarrow H^1 ( G_2, \E[p^2] ) \rightarrow H^1 ( \langle \gamma \rangle, \E[p^2] ).
\end{equation}
But $H^1_{{\rm loc}} ( G_2, \E[p^2] )$ is the intersection of the kernels of the restriction maps $H^1 ( G_2, \E[p^2] ) \rightarrow H^1 ( C^\prime, \E[p^2] )$, as $C^\prime$ varies over all cyclic subgroups of $G_2$ (see  Definition \ref{defloc}).
Since $H^1_{{\rm loc}} ( G_2, \E[p^2] ) \neq 0$, \begin{equation}\label{eqn:relspeciale}
H^1 ( \Gal ( K_1 ( \zeta_{p^2} ) / k ), \E[p^2]^{\langle \gamma \rangle} ) \neq 0.
\end{equation}

\noindent We shall show that  one of the points in the set $\{ P_1, P_2 \}$ is $k$-rational.
 By relation~(\ref{eqn:rel21}), the group $\Gal ( K_1 ( \zeta_{p^2} ) / k )$ is cyclic.
Let $\delta$ be its generator and denote by $\overline{\delta}$ its restriction to $K_1$.
Then
\begin{equation}\label{eqn:rel23}
\overline{\delta} =
\left(
\begin{array}{cc}
\lambda_1 & 0 \\
0 & \lambda_2 \\
\end{array}
\right)
\end{equation}
 generates $G_1$ by (\ref{eqn:rel21}). We are going to show that either $\lambda_1 = 1$ or $\lambda_2 = 1$. Then either $P_1$ or $P_2$ is $k$-rational of exact order $p$.
Suppose  $\lambda_1\neq 1$ and $\lambda_2\neq 1$. We shall contradict relation~(\ref{eqn:relspeciale}). To this end, we compute $H^1 ( \Gal ( K_1 ( \zeta_{p^2} ) / k ), \E[p^2]^{ \langle \gamma \rangle } )$.
Since $\Gal ( K_1 ( \zeta_{p^2} ) / k )$ is cyclic and generated by $\delta$, then Lemma~\ref{lem:lem21} gives
\begin{equation}\label{eqn:rel24}
H^1 ( \Gal ( K_1 ( \zeta_{p^2} ) / k ), \E[p^2]^{ \langle \gamma \rangle } ) \cong \ker ( N_{\Gal ( K_1 ( \zeta_{p^2} ) / k )} ) / {\rm Im} ( \delta - 1 ).
\end{equation}

\noindent We first show that
\begin{equation}\label{eqn:rel25}
{\rm Im} ( \delta - 1 ) = \E[p^2]^{ \langle \gamma \rangle }.
\end{equation}

\noindent Let $\widetilde{\delta}$ be an extension of $\delta$ to $K_2$.
Then, by~(\ref{eqn:rel23}), there exist $\mu_1, \mu_2, \mu_3, \mu_4 \in \F_p$ such that
\[
\widetilde{\delta}
=
\left(
\begin{array}{cc}
\lambda_1 & 0 \\
0 & \lambda_2 \\
\end{array}
\right)
+
p
\left(
\begin{array}{cc}
\mu_1 & \mu_2 \\
\mu_3 & \mu_4 \\
\end{array}
\right).
\]
Since $\lambda_1 \not \equiv 1 \mod ( p )$ and $\lambda_2 \not \equiv 1 \mod ( p )$, $\widetilde{\delta} - I$ is invertible.
Thus
\[
\widetilde{\delta} - I \colon \E[p^2] \rightarrow \E[p^2]
\]
is an isomorphism.
Since, $\delta$ is a restriction of $\widetilde{\delta}$ to $K_1 ( \zeta_{p^2} )$, then \[
\delta - I \colon \E[p^2]^{\langle \gamma \rangle} \rightarrow \E[p^2]^{\langle \gamma \rangle}
\]
is an isomorphism.
In particular, $\delta - I$ is surjective, which gives relation~(\ref{eqn:rel25}).
By~(\ref{eqn:rel24}), we conclude
\[
H^1 ( \Gal ( K_1 ( \zeta_{p^2} ) / k ), \E[p^2]^{ \langle \gamma \rangle } ) = 0,
\]
contradicting relation~(\ref{eqn:relspeciale}).

Finally, the elements of  $G_1$ which fix $\zeta_p$ are diagonal and have determinant $1$. In addition, we just proved that every element of $G_1 $ fixes $P_1$ or $P_2$. But only the identity has these three properties, thus $k ( \zeta_p ) = K_1$. \CVD

\begin{pro}[The non-cyclic case]\label{pro:pro31}
Suppose that $k$  does not contain the subfield  $L$ of $\Q ( \zeta_{p^2} )$ such that $[L : \Q] = p$.
Suppose that  $H^1_{{\rm loc}} ( G_2, \E[p^2] ) \neq 0$ and that $G_1$ is not cyclic.
Then, there exists a $k$-rational point of $\E$ of exact order $p$.
Moreover $K_1 = k^\prime ( \zeta_p )$, where $k^\prime$ is the subfield of $K_1$ defined in Lemma~\ref{lem:lem14}(ii).
\end{pro}

\DIM By Lemma~\ref{lem:lem14}(ii), there exists a $\F_p$-basis $\{P_1, P_2 \}$ of $\E[p]$ such that  $G_1$ is generated by \[
\rho =
\left(
\begin{array}{cc}
\lambda_1 & 0 \\
0 & \lambda_2 \\
\end{array}
\right)
\ {\rm and }\
\sigma =
\left(
\begin{array}{cc}
1 & 1 \\
0 & 1 \\
\end{array}
\right).
\]

We shall prove that  $P_1$ is defined over $k$, or equivalently that $\lambda_1=1$.  Suppose $\lambda_1\not=1$, then in particular $\rho$ is not the identity and its order  $d$ is greater than $1$. Recall  that
 $d$ divides $p-1$, while the oder of $\sigma $ is $p$. Then $[K_1:k]=dp$.\\

 Consider the  field $F \subseteq K_1$  fixed by $\langle \sigma \rangle$. The point $P_1$ is fixed by $ \sigma$ and therefore it is defined over $F$.
An easy computation shows that $\langle \sigma \rangle$ is normal in $G_1$. Then, the  field $F $ is normal over $k$. In addition  $\Gal ( F/ k )$ is cyclic, because it is generated by a restriction of $\rho$ to $F$.
Since $d$ and $p$ are coprime,  the restriction of $\rho$ to $F$ has still order $d$.
Hence
\begin{equation}\label{eqn:rel31}
\Gal ( F/ k ) = \langle \rho \rangle\cong \Z / d \Z \ \hspace{0.2cm} {\rm and} \hspace{0.2cm} \ 1\not=d \ {\rm divides}  \ p-1.
\end{equation}

Consider  $F ( \zeta_{p^2} ) / k$.
Since $k$ does not contain $L$, $[kL : k] = p$.  By relation~(\ref{eqn:rel31}), $[F: k] = d$, prime to $p$. Therefore $[F L : k] = d p$.
As $\sigma$ has determinant $1$,  then $ \zeta_p  \in F$.
Hence  $F L = F ( \zeta_{p^2} )$.
The extensions $F/ k$ and $kL / k$ are Galois and $k$-linearly disjoint. Hence
\begin{equation}\label{eqn:rel32}
\Gal ( F ( \zeta_{p^2} ) / k ) \cong \Gal ( kL / k ) \times \Gal ( F / k ) \cong \Z / pd \Z=\langle \delta\rangle.
\end{equation}

Let  $\overline{\delta}$ be the  restriction of $\delta$ to $F$. Since  ${\delta}$ generates $\Gal ( F ( \zeta_{p^2} ) / k )$, then $\overline{\delta}$ generates $\Gal ( F / k )$. Then  $\overline{\delta}= \rho^t$, with $t$ prime to $d$. Thus
\begin{equation}\label{eqn:rel36}\delta ( P_1 ) = \overline{\delta} ( P_1 ) = \lambda P_1,   \ { \rm with }\  1\not=\lambda=\lambda^t_1 \in \F_p^\ast.\end{equation}

We shall show that the hypothesis  $\lambda\not=1$ implies that  a cohomological group is both trivial and nontrivial, which is absurd.

\noindent Consider the group $H^\prime = \Gal ( K_2 / F ( \zeta_{p^2} ) )$. By relation (\ref{eqn:rel32})
\[
[F ( \zeta_{p^2} ) : k] = d p = [K_1: k].
\]
Therefore $\vert H^\prime \vert=\vert H \vert$, with $H = \Gal ( K_2 / K_1 )$ of cardinality $p^2$, by Lemma~\ref{lem:lem16}. Hence $\vert H^\prime \vert = p^2$ and $H'$ is a $p$-group.
In addition $H^\prime$ is normal in $G_2$, because $F ( \zeta_{p^2} ) / k$ is a Galois extension.

Consider the following inflaction-restriction sequence:
\begin{equation*}\label{eqn:rel33}
0 \rightarrow H^1 ( \Gal ( F ( \zeta_{p^2} ) / k ), \E[p^2]^{H^\prime} ) \rightarrow H^1 ( G_2, \E[p^2] ) \rightarrow H^1 ( H^\prime, \E[p^2] ).
\end{equation*}

 We first show that $H^1_{{\rm loc}} ( H^\prime, \E[p^2] ) = 0$.
Let $\pi \colon \textrm{GL}_2 ( \Z / p^2 \Z ) {\rightarrow} \textrm{ GL}_2 ( \Z / p \Z )$
be the reduction modulo $p$.
Then $\ker ( \pi ) \cap H^\prime=H \cap H^\prime$ is a $\F_p$-vector space and its dimension is $1$. Indeed
 $H \cap H^\prime = \Gal ( K_2 / K_1 ( \zeta_{p^2} ) )$, which is a cyclic group of order $p$.
Thus, the intersection between the $p$-Sylow subgroup of $H^\prime$ (which is   $H^\prime$ itself) and $\ker ( \pi )$, has dimension $1$ over $\F_p$.
By~\cite[Proposition 3.2]{DZ}, we have
\begin{equation*}\label{eqn:rel34}
H^1_{{\rm loc}} ( H^\prime, \E[p^2] ) = 0.
\end{equation*}

\noindent  We\hspace{0.1cm} remark\hspace{0.1cm} that\hspace{0.1cm}  the\hspace{0.1cm} restriction\hspace{0.1cm} $H^1 ( G_2, \E[p^2] ) \rightarrow H^1 ( H^\prime, \E[p^2] )$\hspace{0.1cm} sends $H^1_{{\rm loc}} ( G_2, \E[p^2] )$ to   $H^1_{{\rm loc}}  ( H^\prime, \E[p^2] )$.
In fact,  if a cocycle satisfies the local conditions relative to $G_2$, it satisfies them relative to any of its subgroups (see  Definition \ref{defloc}).
By hypothesis  $H^1_{{\rm loc}} ( G_2, \E[p^2] ) \neq 0$ and we just proved that $H^1_{{\rm loc}} ( H^\prime, \E[p^2] ) = 0$.
Then, by the exact sequence,
\begin{equation}\label{eqn:rel35}
H^1 ( \Gal ( F ( \zeta_{p^2} ) / k ), \E[p^2]^{H^\prime} ) \neq 0.
\end{equation}

On the other hand,  by Lemma~\ref{lem:lem21},
\begin{equation*}\label{eqn:rel37}
H^1 ( \Gal ( F ( \zeta_{p^2} ) / k ), \E[p^2]^{ H^\prime } ) \cong \ker ( N_{\Gal ( F ( \zeta_{p^2} ) / k )} ) / {\rm Im} ( \delta - 1 ).
\end{equation*}

\noindent We are going to show that
\begin{equation*}\label{eqn:rel38}
{\rm Im} ( \delta - 1 ) = \E[p^2]^{ H^\prime }.
\end{equation*}
This implies \[
H^1 ( \Gal ( F ( \zeta_{p^2} ) / k ), \E[p^2]^{ H^\prime } ) = 0
\]
and contradicts (\ref{eqn:rel35}).
 The group $H^\prime$ fixes the field $F ( \zeta_{p^2} )$. So  $P_1 \in \E[p^2]^{H^\prime}$,  because $P_1$ is defined over $ F$. On the other hand  $F ( \zeta_{p^2} )$ does not contain $K_1$, indeed  $\Gal ( F ( \zeta_{p^2} ) / k )$ is cyclic of order $dp$ and $G_1$  has the same order, but $G_1$ is not cyclic because $d>1$. Therefore $\E[p^2]^{H^\prime}$ does not contain $\E[p]$ and
either $\E[p^2]^{H^\prime} = \langle P_1 \rangle, \ {\rm or} \ \E[p^2]^{H^\prime} =\langle Q_1 \rangle$, where $p Q_1 = P_1$.
In any case  $\E[p^2]^{H^\prime} = \langle a Q_1 \rangle$, with $a$  a non-zero element of $\Z/ p^2 \Z$.
By~(\ref{eqn:rel36}), there exists $\mu \in \F_p$ such that
\[
\delta - I ( a Q_1 ) = ( \lambda + p \mu - 1 ) a Q_1.
\]
Since $\lambda \not \equiv 1 \mod ( p )$, then $( \lambda + p \mu - 1 ) a Q_1$ generates $\E[p^2]^{H^\prime}$.
Therefore
\[
\delta - I \colon \E[p^2]^{H^\prime} \rightarrow \E[p^2]^{H^\prime}
\]
is surjective and ${\rm Im} ( \delta - 1 ) = \E[p^2]^{H^\prime}$, as desired.\\

In conclusion $\rho= {1 \ 0\choose 0\  \lambda_2}$. Let us show that $K_1 = k^\prime ( \zeta_p )$, with  $k^\prime $  the field fixed by $\rho$. If $\lambda_2=1$,  then $K_1=k'$, but $\zeta_p \in K_1$. If $\lambda_2\not=1$, then  $\zeta_p \not\in k'$, indeed only elements of determinant $1$ fix $\zeta_p$. Thus $  K_1 \supset k^\prime(\zeta_p )$. In addition the two extensions have the same degree, so they are equal. \CVD

\noindent {\bf Proof of Theorem~\ref{main}.} By Propositions~\ref{pro:pro22} and~\ref{pro:pro31},
we have that if $\E$ does not have any $k$-rational point of order $p$, then $H^1_{{\rm loc}} ( G_2, \E[p^2] ) = 0$.
By Theorem \ref{si1}, if $H^1_{{\rm loc}} ( G_2, \E[p^2] ) = 0$, then every point $P \in \E ( k )$, which is locally divisible by $p^2$ in $\E ( k_v )$ for all but finitely many primes $v$, is divisible by $p^2$ in $\E ( k )$.
\CVD

\newpage

\thispagestyle{empty}

Laura Paladino\par\smallskip
Dipartimento di Matematica \par
Universit\`{a} della Calabria\par\smallskip
via Ponte Pietro Bucci, cubo 31b\par
IT-87036 Arcavacata di Rende (CS)\par
e-mail address: paladino@mat.unical.it

\vskip 1cm

Gabriele Ranieri\par\smallskip
Mathematisches institut der\par
Universit\"{a}t G\"{o}ttingen\par\smallskip
Bunsenstrasse 3-5,\par
D-37038 G\"{o}ttingen\par
e-mail address: ranieri@mat.unicaen.fr

\vskip 1cm

Evelina Viada\par\smallskip
Departement Mathematik\par
Universit\"{a}t Basel\par\smallskip
Rheinsprung, 21\par
CH-4051 Basel\par
e-mail address: evelina.viada@unibas.ch

\end{document}